\documentclass[11pt,a4paper]{article}
\usepackage{amssymb}
\usepackage{subfloat,amsmath,graphics,amssymb,amsfonts,graphicx,lineno}
\usepackage{amsmath}
\usepackage{amssymb}
\usepackage{amsthm}
\usepackage{amsfonts}
\newtheorem{theorem}{Theorem}[section]

\newtheorem{lemma}[theorem]{Lemma}
\newtheorem{proposition}[theorem]{Proposition}
\theoremstyle{definition}

\theoremstyle{definition}

\begin{document}
\begin{center}
\textbf{Notes on amenability}
\\ \ \\
Miad Makareh Shireh
\ \\ \ \\
\end{center}
 ABSTRACT: We show that for a Banach algebra $A$ with a bounded approximate identity, the amenability of $A\widehat{\otimes}A$, the amenability of $A\widehat{\otimes}A^{op}$ and the amenability of $A$ are equivalent. Also if $A$ is a closed ideal in a commutative Banach algebra $B$, then the weak amenability of $A\widehat{\otimes}B$ implies the weak amenability of $A$.
 \section{Introductions and Preliminaries}
 Let $A$ to be a Banach algebra and $X$ an $A$-bimodule that is a Banach space. We say that $X$ is a Banach $A$-bimodule if there exists  constant $C>0$ such that
 \begin{align*}
 \|a.x\|&\leq C\|a\|\|x\|,\\
 \|x.a\|&\leq C\|a\|\|x\| \hskip 1 cm (a \in A, x \in X).
 \end{align*}
 If $X$ is a Banach $A$-bimodule, then $X^*$ is a Banach $A$-bimodule for the actions defined by \begin{align*}
\langle a.f,x\rangle &=\langle f,x.a\rangle\\ \langle f.a,x\rangle &=\langle f,a.x\rangle \hskip 0.5cm (a\in A, f\in X^*,x \in X). \end{align*}
The Banach $A$-bimodule $X^*$ defined in this way is said to be a dual Banach $A$-bimodule.\\\\

 A linear mapping $D$ from $A$ into $X$ is a derivation if $$ D(ab)=a.D(b)+D(a).b \hskip 1 cm (a,b \in A).$$
For $x \in X$, the mapping $ad_x:A\longrightarrow X$ defined by $ad_x(a)=a.x-x.a$ is a continuous derivation. The derivation $D$ is inner if there exists $x \in X$ such that $D=ad_x$.\\
 $A$ is said to be amenable if for every Banach $A$-bimodule $X$ , any continuous derivation from $A$ into the dual Banach $A$-bimodule $X^*$ is inner. This notion has been introduced in [4] and has been studied extensively since.\\
 The Banach algebra $A$ is said to be weakly amenable if  any continuous derivation from $A$ into the dual Banach $A$-bimodule $A^*$ is inner. This notion was first introduced in [1] for the commutative case and then in [6] for the general case.
 Every Amenable Banach algebra has a bounded approximate identity [ 8, Proposition 2.21]. Also if $A$ and $B$ are two amenable Banach algebras then also is $A\hat{\otimes}B$ [ 4, Proposition 5.4]. \\
 The question is whether the converse is true or not. The only work on this question is done by B.E. Johnson in the following:
 \begin{proposition}
 Suppose that $A$ is a Banach algebra and $B$ is another Banach algebra such that there exists $b_0 \in B$ with $b_0 \notin \overline{\rm{Lin}\{bb_0-b_0b: b \in B\}}$. If $A\widehat{\otimes}B$ is amenable then $A$ is amenable.
 \end{proposition}
 {\bf Proof:}  \hskip 0.2 cm See [5, Proposition 3.5] \hskip 5 cm $\Box$\\\\
 But still the question remains open for general $A$ and $B$ even for the case $A=B$. \\\\
 In section $2$, we prove that the amenability of $A\hat{\otimes}A$ implies amenability of $A$ in the case that $A$ has a bounded approximate identity. Indeed we show that for a Banach algebra $A$ with a bounded approximate identity the following are equivalent: \begin{itemize}
                 \item [(i)] $A$ is amenable;
                 \item [(ii)] $A\widehat{\otimes}A$ is amenable;
                 \item [(iii)] $A\widehat{\otimes}A^{op}$ is amenable (Where as usual $A^{op}$ is the Banach algebra obtaining by reversing the product of $A$).
               \end{itemize}
               Since having a bounded approximate identity is a necessary condition for amenability, we can not omit the condition that $A$ has  bounded approximate identity unless we can prove that amenability of  $A\hat{\otimes}A$ necessitates having a bounded approximate identity for $A$.
               \\\\
 In section 3 we investigate the question  but for weak amenability instead of amenability. We prove that if $B$ is a commutative Banach algebra and $A$ is a closed ideal in $B$ , then the weak amenability of  $A\hat{\otimes}B$ implies the weak amenability of $A$.
 \section{The amenability results}
In this section we try to answer the question whether amenability $A\hat{\otimes}B$ implies the amenability of $A$ and $B$ or not. We mainly concentrate  on the special case where $A=B$. However, we will also obtain some results about the case where $A$ is not necessarily equal to $B$.

\noindent First we start with a simple result:
\begin{theorem}
Suppose that $A$ and $B$ are Banach algebras and $B$ has a non- zero character. If $A\widehat{\otimes}B$ is amenable, then $A$ is also amenable.
\end{theorem}
{\bf Proof:} Let $\varphi$ be a non- zero character in $B$ and define the unique mapping $\theta:A\widehat{\otimes}B \longrightarrow A$ acting on elementary tensors by $$\theta(a\otimes b)=\varphi(b)a \hskip 0.5 cm (a \in A , b \in B).$$
We show that $\theta$ is an algebra homomorphism (obviously $\theta$ is continuous). Since $\theta$ is linear, it is enough to check this for elementary tensors. To see this we have
 $$\theta((a\otimes b)(c\otimes d))=\theta((ac\otimes bd)=\varphi(bd)ac.$$
On the other hand
 $$ \theta((a\otimes b))\theta((c\otimes d))=\varphi(b)a\varphi(d)c=\varphi(bd)ac.$$
So
 $$\theta((a\otimes b)(c\otimes d))=\theta((a\otimes b))\theta((c\otimes d)).$$
 And since $\varphi$ is non - zero,  $\theta$ is surjective and hence  $A$ is amenable. \hskip 2 cm  $\Box$\\\\
Throughout the following we let $ \pi: A\widehat{\otimes}A^{op}\longrightarrow A$ be the so-called product map; mapping specified by acting on elementary tensors by $\pi(a\otimes b)=ab \hskip 0.5 cm (a,b \in A)$ and we let $K={\rm ker } \pi$.\\
The Banach algebra $A$ can be made into a left $ A\widehat{\otimes}A^{op}$-module by the module multiplication specified by  $$(a\otimes b).c= acb \hskip 1 cm (a,b,c \in A).$$
\begin{theorem}\label{1}
Suppose that $A\widehat{\otimes}A^{op}$ is amenable and $A$ has a
bounded approximate identity. Then $A$ is amenable.
\end{theorem}
{\bf Proof:} Since $A$ has a bounded approximate identity, the
short exact sequence \hskip 0.5cm $(\prod^{op})^*:0 \longrightarrow
A^* \stackrel{\pi^*} \longrightarrow (A\widehat{\otimes}A^{op})^*
\stackrel{\imath^*}\longrightarrow K^*\longrightarrow 0$ is an
admissible short exact sequence of right
$A\widehat{\otimes}A^{op}$-modules.
($\imath$ is the inclusion map).\\
Since $A^*$ is a dual $A\widehat{\otimes}A^{op}$-module, from [2, Theorem 2.3] , $(\prod^{op})^*$ splits and since
$A\widehat{\otimes}A^{op}$ has a bounded approximate identity and $\pi$ is onto
,  [2, Theorem 3.5]
implies that $K$ has a bounded right approximate identity. Now since
$A$ has a bounded approximate identity, from  [2, Theorem 3,10]
$A$ is amenable. \hskip 1 cm $\Box$ \\\\
Theorem 2.3 has been the motivation for us to consider the question of under which conditions on the tensor products, $A$ has a
bounded approximate identity. The following is one of them. Before going to next Theorem, we need a Lemma.
\begin{lemma}
Let $A$ to be a Banach algebra with a two-sided bounded approximate identity and $X$ a Banach $A$-bimodule on which $A$ acts trivially on one side. Then for every continuous derivation $D$ from $A$ into $X$, there exists a bounded net $(\zeta_i)_i$ in $X$ such that $D(a)=\lim_i a.\zeta_i-\zeta_i.a \hskip 0.5 cm  (a \in A)$.
\end{lemma}
{\bf Proof: } Since we can embed $X$ into $X^{**}$ through the canonical injection, we can consider $D$ as a continuous derivation into the dual module $X^{**}$. Also since the action of $A$ on one side of $X$ is trivial, action of $A$ on other side of $X^*$ is trivial. Therefore $D$ is inner. Hence there exists $\xi \in X^{**}$ such that $$D(a)=a.\xi-\xi.a \hskip 1 cm (a \in A).$$
Now by Goldstein's Theorem, there is a bounded net $(\tau_j)_{j\in J}$ in $X$ converging to $\xi$ in weak$^*$ topology of $X^{**}$.
Thus $$D(a)=a.\xi-\xi.a=\rm{wk}^*-\lim_j a.\tau_j-\tau_j.a \hskip 1 cm  (a \in A),$$
and hence $$D(a)=\rm{wk}-\lim_j a.\tau_j-\tau_j.a \hskip 1 cm (a \in A). $$
Let $\Delta=\{a_1,a_2,...,a_n\}$ be a finite subset of $A$. Then in $\bigoplus_{i=1}^n X$,  we have $$(D(a_1),...,D(a_n)) \in weak-cl(co(\{(a_1.\tau_j-\tau_j.a_1,...,a_n.\tau_j-\tau_j.a_n): j \in J\}))$$
Therefore by Mazur's Theorem $$ (D(a_1),...,D(a_n)) \in \rm{norm}-\rm{cl(co}(\{(a_1.\tau_j-\tau_j.a_1,...,a_n.\tau_j-\tau_j.a_n) ): j \in J\})$$
And hence for $\epsilon>0$, there exists $\zeta_{\Delta,\epsilon} \in\rm{co}(\{\tau_j:h\in J\})$ such that $$\|D(a_i)-(a_i.\zeta_{\Delta,\epsilon}-\zeta_{\Delta,\epsilon}.a_i)\|<\epsilon \hskip 1 cm  (a_i\in \Delta)$$
So by ordering the set of the finite subsets of $A$ by inclusion and positive real numbers by decreasing order , the net $(\zeta_{\Delta,\epsilon})$ is the desired net. \hskip 1 cm $\Box$
\begin{theorem}\label{2}
Suppose that $A\widehat{\otimes}A^{op}$ has a bounded approximate
identity and each one of the topologies on $A$ defined by the family of
seminorms $\rho_a: b\mapsto\|ab\|$ and $\gamma_a:b\mapsto\|ba
\|$ is stronger than weak topology on $A$. Then $A$ has a
(two-sided) bound approximate identity.
\end{theorem}
{\bf Proof}: Suppose that $A\widehat{\otimes}A^{op}$ has a
bounded approximate identity. we consider $A$ as an
$A\widehat{\otimes}A^{op}$-bimodule by actions specified by:
\begin{align*}
(a\otimes b)\bullet c &=acb \\
c\bullet (a\otimes b) &=0 \hskip 0.5cm (a,b,c \in A)
\end{align*}
 It can be easily seen that $A$ is
a Banach $A\widehat{\otimes}A^{op}$-bimodule by the actions above . Now we define a
derivation $D:A\widehat{\otimes}A^{op}\longrightarrow A$ by acting on elementary tensors as $D(a\otimes b)= ab \hskip 0.5 cm (a,b\in A)$. $D $ is obviously continuous and also $D$ is a
 derivation since
 $$D((a\otimes b)\cdot (c\otimes d))= D(ac\otimes db)= acdb$$
 ($\cdot$ is the product in $A\widehat{\otimes}A^{op}$ ). On the other hand:
 $$(a\otimes b)\bullet D(c\otimes d)+ D(a\otimes b)\bullet (c\otimes d)= (a\otimes
 b)\bullet cd = acdb$$
 Therefore $D \in \mathbb{Z}^1(A\widehat{\otimes}A^{op},A)$.
 Now since the right action of $A\widehat{\otimes}A^{op}$ on
 $A$ is trivial and $A\widehat{\otimes}A^{op}$ has a bounded approximate identity, from Lemma 2.4,  there exists
 a bounded net$(\zeta_i)_i$ in $A$ such that
 $D(a\otimes b)=\lim_i
 ad_{\zeta_i}(a\otimes b)$.\\
 So $ab= \lim_i a\zeta_ib \hskip 0.5 cm (a, b \in A)$. Thus for all $a,b \in A$
 \begin{equation}
 \lim_i a(b-\zeta_ib)=0 \hskip 1cm
 \lim_i(b-b\zeta_i)a=0
\end{equation}
If we denote the topology induced by the family of seminorms
$\{\rho_a| a\in A\}$ by $\tau$ and the topology induced by the
family of seminorms $\{ \gamma_a|a\in A\}$ by $\varsigma$, then from (1) we have\\
$$a\zeta_i\longrightarrow a    \hskip 0.5cm ( in \hskip 0.1cm
\tau \hskip 0.1cm for\hskip 0.2cm all\hskip 0.1cm   a \in A) \hskip 3cm (2)$$\\\\
$$\zeta_ia\longrightarrow a \hskip 0.5cm ( in \hskip 0.1cm
\varsigma \hskip 0.1cm for\hskip 0.2cm all\hskip 0.1cm   a \in
A) \hskip 3cm (3)$$\\\\
 since we assume both $\tau$ and $\varsigma$ to be stronger than weak
topology on $A$, then by (2) and (3), $A$ has a weakly two-sided
bounded approximate identity and hence $A$ has a two-sided
bounded approximate identity.\hskip 6 cm $\Box$\\\\
\begin{theorem}\label{3}
Suppose that $A\widehat{\otimes}A^{op}$ is amenable and that $A$
has the property that each one of the topologies induced on $A$ by the family of
seminorms $\{\rho_a| a\in A\}$ where $\rho_a(b)=\|ab\|$ and $\{\gamma_a| a\in A\}$ where
$\gamma_a(b)=\|ba\| $, are stronger than the weak
topology on $A$. Then $A$ is amenable.
\end{theorem}
{\bf Proof:} Firstly by the fact that
$A\widehat{\otimes}A^{op}$ necessarily has a two-sided bounded
approximate identity and from Theorem 2.5 we have that $A$ has a two-sided bounded
approximate identity and then from Theorem 2.3 we have $A$ is
amenable.\hskip 1 cm $\Box$\\\\
In next Theorem we attempt to relate amenability of
$A\widehat{\otimes}A$ (in the case that $A$ has a bounded approximate identity) to
the amenability of $A\widehat{\otimes}A^{op}$ and then by using the
preceding theorems, we attempt to prove the amenability of $A$
when $A\widehat{\otimes}A$ is amenable.
Before going to next Theorem, we need a Lemma.
\begin{lemma}
Let $A$  be Banach algebra with a bounded approximate identity such that for any neo-unital Banach $A$-bimodule $X$ and $Y$ a closed submodule of $X$, every $f\in Z_A(Y^*)$ can be extended to a functional $\tilde{f} \in Z_A(X^*)$. Then $A$ is amenable.
\end{lemma}
{\bf proof:} As in the proof of [7, Theorem 1], for concluding the amenability of
$A$, it is enough to have the property in the Lemma for the Banach $A$-bimodule  $L=(A\widehat{\otimes}A)^*\widehat{\otimes} (A\widehat{\otimes}A)$ with the module actions specified by
\begin{align*}
a.(x^*\otimes x)&=x^*\otimes a.x,
\\ (x^*\otimes x).a&=x^*\otimes x.a \hskip 0.5 cm (a \in A, x \in (A\widehat{\otimes}A), x^*\in (A\widehat{\otimes}A)^*).
\end{align*}
Since $A$ has bounded approximate identity, $X=A\widehat{\otimes}A$ is neo-unital and hence by the above definition of the actions of $A$ on $L$, $L$ is also neo-unital. \hskip 1 cm $\Box$

\begin{theorem} \label{5}
Suppose that $A$ is a Banach algebra with a bounded approximate
identity such that $A\widehat{\otimes}A$ is amenable. Then
$A\widehat{\otimes}A^{op} $ is also amenable.
\end{theorem}
{\bf proof:} Suppose that $X$ is a Banach neo-unital
$A\widehat{\otimes}A^{op}$-bimodule and that $\bullet$ denotes the action of $A\widehat{\otimes}A^{op}$ on $X$.  We define:
\begin{align*}
 (a\otimes b)\circ x &= \lim_i(a\otimes
e_i)\bullet x \bullet(e_i\otimes b),\\
 x \circ (a\otimes b) &= \lim_i(e_i \otimes b)\bullet x \bullet( a\otimes e_i)\hskip
 0.5cm (x\in X  \rm{and} \hskip 0.2cm  a,b\in A).
 \end{align*}
 First we note  that the above limits exist because by the
 assumption that $X$ is neo-unital we have: \\
 If $x \in X$ then there exist $y \in X$ and $ u,v \in
 A\widehat{\bigotimes}A^{op}$ such that $x=u\bullet y\bullet v$ and then we
 have:
$$(a\otimes e_i)\bullet x \bullet(e_i\otimes b)=(a\otimes
e_i)\bullet u\bullet y\bullet v \bullet(e_i\otimes b)=((a\otimes
e_i)\star u)\bullet y\bullet (v\star(e_i\otimes b)),$$ where $\star$ denotes the product in $A\widehat{\otimes}A^{op}$. Since
$(e_i)_{i\in \Lambda}$ is a bounded approximate identity for $A$, it
can be easily seen that $\lim_i (a\otimes e_i)\star u = a.u$ and
$\lim_i v\star (e_i\otimes b)= v.b,$ where $a.(e\otimes f)= ae\otimes
f $ and $(e\otimes f).b=e\otimes bf$.\\
So $ \lim_i(a\otimes e_i)\bullet x \bullet(e_i\otimes b)$ exists and
we can similarly prove the existence of the second limit.
Also $\circ$ induces a module action of $A\widehat{\bigotimes}A$  on $X$ .
To see the reason, by linearity, it is enough to check the module conditions for elementary tensors.
$$((a\otimes b)(c\otimes d))\circ x = (ac\otimes bd)\circ x =
\lim_i(ac\otimes e_i)\bullet x \bullet
 (e_i\otimes bd)$$

\noindent On the other hand:
\begin{align*}
 (a\otimes b)\circ ((c\otimes d)\circ x) &= (a\otimes b)\circ
 (\lim_j(c\otimes e_j)\bullet x \bullet  (e_j\otimes d))\\
&=\lim_i(a\otimes e_i)\bullet(\lim_j(c\otimes e_j)\bullet x \bullet (e_j\otimes d))\bullet(e_i\otimes
 b)\\ &= {\rm lim_ilim_j} (ac\otimes e_je_i)\bullet x \bullet(e_je_i\otimes bd)\\
 &= {\rm lim_i} (ac\otimes e_i)\bullet x \bullet(e_i\otimes bd).
 \end{align*}
Hence $$((a\otimes b)(c\otimes d))\circ x = (ac\otimes bd)\circ x= (a\otimes b)\circ ((c\otimes d)\circ x).$$
In a similar way we can show that $$ x \circ ((a\otimes b)(c\otimes d))=(x\circ(a\otimes b))(c\otimes d).$$
Also we have:
\begin{align*}
((a\otimes b)\circ x)\circ(c\otimes d) &= {\rm lim_i} (e_i\otimes
 d)\bullet({\rm lim_j} (a\otimes e_j)\bullet x\bullet(e_j\otimes b))\bullet (c\otimes e_i)\\
&={\rm lim_ilim_j} ((e_i\otimes d)\star(a\otimes e_j)) \bullet x \bullet((e_j\otimes
 b)\star(c\otimes e_i))\\ &={\rm lim_ilim_j} (e_ia\otimes e_jd) \bullet x\bullet
(e_jc\otimes e_ib)\\ &=(a\otimes d)\bullet x\bullet (c\otimes b).
\end{align*}
On the other hand:
\begin{align*}
(a\otimes b)\circ(x \circ(c\otimes d)) &={\rm lim_ilim_j} (a\otimes
e_i)\bullet((e_j\otimes d)\bullet x\bullet (c\otimes e_j))\bullet
(e_i\otimes b)\\ &={\rm lim_ilim_j}((ae_j\otimes de_i)\bullet x\bullet (ce_i\otimes be_j)\\
&=(a\otimes d)\bullet x\bullet (c\otimes b).
\end{align*}
Hence $$((a\otimes b)\circ x)\circ(c\otimes d)=(a\otimes b)\circ(x \circ(c\otimes d)).$$
 So $X$ is an $A\widehat{\otimes}A$-bimodule for the action $\circ$. Also since the net $(e_i)$ is bounded, it can be easily seen that $X$ is indeed a Banach $A\widehat{\otimes}A$-bimodule for $\circ$. For a Banach
 $A\widehat{\otimes}A^{op}$-bimodule $X$,  $X_\dag$ denotes $X$ as an  $A\widehat{\otimes}A$-bimodule (via the action $\circ$).\\
 Now if $Y$ is a closed submodule of
$X$ and $f \in
 Z_{A\widehat{\otimes}A^{op}}(Y^*)$, we show that \\
 $f \in Z_{A\widehat{\otimes}A}(Y_\dag^*)$.\\ To prove the above
statement we have:
\begin{align*}
(a\otimes b)\circ f &=\rm{wk}^*-\lim_i(a\otimes e_i)\bullet f \bullet
(e_i\otimes b)\\ &= \rm{wk}^*-\lim_if\bullet (a\otimes e_i)\bullet(e_i\otimes
b)\\ &=\rm{wk}^*-\lim_if\bullet(ae_i\otimes be_i)\\ &= f\bullet (a\otimes b).
\end{align*}
 Similarly $$f \circ (a\otimes b)=(a\otimes b)\bullet f.$$
 Thus
   $$f \in Z_{A\widehat{\otimes}A}(Y_\dag^*).$$

\noindent Now from [7, Theorem 1] , $f$ has an extension to an
$\tilde{f} \in Z(A\widehat{\otimes}A , X_\dag^*)$ .\\
 We show that $\tilde{f}\in Z(A\widehat{\otimes}A^{op}, X^*)$
For this purpose we have
\begin{align*}
(a\otimes b)\bullet \tilde{f} &=\rm{wk}^*-\lim_i-\rm{wk}^*-\lim_j ((a\otimes
e_i)(e_j\otimes b))\bullet \tilde{f} \bullet (e_i\otimes
e_j)\\ &= \rm{wk}^*-\lim_i(a\otimes e_i) (\rm{wk}^*-\lim_j(e _j\otimes
b)\bullet \tilde{f} \bullet (e_i\otimes e_j))\\
 &= \rm{wk}^*-\lim_i (a\otimes e_i)\bullet(\tilde{f}\circ (e_i\otimes b))\\ &=\rm{wk}^*-\lim_i
(a\otimes e_i)\bullet ((e_i\otimes b)\circ \tilde{f})
\\ &=\rm{wk}^*-\lim_i (a\otimes e_i)\bullet(\rm{wk}^*-\lim_j(e_i\otimes
e_j)\bullet \tilde{f}\bullet (e_j\otimes
b))\\ &=\rm{wk}^*-\lim_i\rm{wk}^*-\lim_j (ae_i\otimes e_je_i)\bullet
\tilde{f} \bullet (e_j\otimes b)\\ &=\rm{wk}^*-\lim_i (a\otimes e_i)\bullet
\tilde{f}\bullet(e_i\otimes b)\\ &=(a\otimes b)\circ \tilde{f}.
\end{align*}
similarly we have $\tilde{f}\bullet(a\otimes b)=\tilde{f}\circ
(a\otimes b)$ and since $\tilde{f} \in
Z_{A\widehat{\otimes}A}(X_\dag^*)$, then $(a\otimes b)\bullet
\tilde{f}=\tilde{f}\bullet (a\otimes b)$. Hence $$\tilde{f}\in
Z_{A\widehat{\otimes}A^{op}}(X^*).$$
Since $Y$ was an arbitrary closed submodule of $X$ and   $f$
was arbitrary in $Z_{A\widehat{\otimes}A^{op}}(Y^*)$, again by
exploiting [7,Theorem 1], we have that $A\widehat{\otimes}A^{op}$ is
amenable. \hskip 1 cm$\Box$

\begin{theorem}\label{6}
Suppose that $A\widehat{\otimes}A$ is amenable and $A$ has a
bounded approximate identity. Then $A$ is amenable.
\end{theorem}
{\bf Proof:} By the preceding Theorem we have that
$A\widehat{\otimes}A^{op}$ is amenable. Since $A$ has a bounded
approximate identity, from  Theorem 2.3 , $A$ is amenable.\hskip 1 cm$\Box$\\

Since having a bounded approximate identity is a necessary condition
for an algebra to be amenable, the Theorem 2.8 has the minimum
conditions. If we can prove that amenability of
$A\widehat{\otimes}A$ implies that $A$ has a bounded approximate
identity, then we can even drop the condition in Theorem 2.8 that $A$
has a bounded approximate identity.
\section{Some results in commutative Banach algebras}
Now we go to the case where our algebra $A$ is commutative. First
we prove the following general result.\\
For the Banach algebra $A$, we define $$A^2=\rm{Lin}\{ab: a,b\in A\}.$$

\begin{theorem} \label{7}
Suppose that $B$ is a Banach algebra and $A$ is a closed subalgebra
of $B$ such that $A\widehat{\otimes} B$ is weakly amenable.
Then $(A^2)^-=A$ \end{theorem} {\bf Proof:} {\it Suppose} that
$(A\widehat{\otimes} B)$ is weakly amenable and $ (A^2)^-\neq A$.
Then from Hahn-Banach Theorem there exists a $\lambda \in A^*$ such
that $\lambda|_{A^2}=0 $ and $\lambda\neq0$. So there exists an
$a_0\in A$ such that $\lambda(a_0)=1$. We denote a Hahn-Banach
extension of $\lambda$ on $B$ by $\tilde{\lambda}$. So
$\tilde{\lambda}\in B^*$ and we specify $D:(A\widehat{\otimes}
B)\longrightarrow (A\widehat{\otimes} B)^*$ by $$D(a\otimes b)=
\tilde{\lambda}(a)\tilde{\lambda}(b)(\tilde{\lambda}\otimes \tilde{\lambda}) \hskip 1 cm (a\in A, b \in B),$$ \\
where $(\tilde{\lambda}\otimes\tilde{ \lambda})(c\otimes d)=
\tilde{\lambda}(c)\tilde{\lambda}(d)$.\\\\
Then we have $$D((a\otimes
b)(c\otimes d))= D(ac\otimes bd)=
\tilde{\lambda}(ac)\tilde{\lambda}(bd)(\tilde{\lambda}\otimes\tilde{ \lambda})=0$$\\
On the other hand for $a,c,x \in A$ and $ b,d,y \in B$ we have
$$\langle(a\otimes b)D(c\otimes d), x\otimes y\rangle=\langle D(c\otimes d),xa\otimes
yb\rangle=\tilde{\lambda}(c)\tilde{\lambda}(d)\tilde{\lambda}(xa)\tilde{\lambda}(yb)=0$$
and similarly $$\langle D(a\otimes b).(c\otimes d), x\otimes
y\rangle=\langle D(a\otimes b),cx\otimes dy
\rangle=\tilde{\lambda}(a)\tilde{\lambda}(b)\tilde{\lambda}(cx)\tilde{\lambda}(dy)=0.$$

So $D: A\widehat{\otimes} B\longrightarrow(A\widehat{\otimes}B)^*$ is a continuous derivation
 and hence from weak amenability of $(A\widehat{\otimes} B)$
it follows that $D=ad(\xi)$ for some $\xi \in (A\widehat{\otimes}
B)^*$.\\\\
So
\begin{align*}
\langle D(a_0\otimes a_0),(a_0\otimes a_0)\rangle &=\langle(a_0\otimes a_0).\xi -
\xi.(a_0\otimes a_0), a_0\otimes a_0\rangle\\
&=\langle\xi ,(a_0^2\otimes a_0^2)-(a_0^2\otimes a_0^2)\rangle\\ &=0
\end{align*}
But we have:
$$\langle D(a_0\otimes a_0),(a_0\otimes
a_0)\rangle=\tilde{\lambda}(a_0)\tilde{\lambda}(a_0)(\tilde{\lambda}\otimes
\tilde{\lambda})(a_0\otimes a_0)=(\tilde{\lambda}(a_0))^4=1.$$ So
we have come up with a contradiction and hence  $(A^2)^-=A$ \hskip
1cm $\Box$
\begin{theorem}
Suppose that $B$ is a commutative Banach algebra and $A$ is an ideal
in $B$ such that $A\widehat{\otimes} B$ is weakly
amenable. Then $A$ is weakly amenable.
\end{theorem}

{\bf Proof} {\it Suppose} that $A\widehat{\otimes} B$ is weakly
amenable. Then we define $\varphi : A\widehat{\otimes} B
\longrightarrow A$ by $\varphi(a\otimes b)=ab$.  It can be
easily seen that $\varphi$ is continuous and is an algebra
homomorphism. Also by Theorem \ref{7} we have $\varphi(A)^-=A$. Hence from [3, Proposition 2.11], $A$ is weakly amenable. \hskip 9 cm $\Box$

  \end{document}